# Max/Min Puzzles in Geometry IV

James M. Parks
Dept. of Math.
SUNY Potsdam
Potsdam, NY
*parksjm@potsdam.edu*
*April 5, 2025*

**Abstract**. In the previous paper, Max/Min Puzzles in Geometry III, we searched for the smallest area triangle which contained a regular unit polygon (Square, Pentagon, Hexagon). In this paper we will work in 3-dimensions, and search for the smallest regular Tetrahedron which contains a regular unit polyhedron (Cube, Octahedron, Icosahedron, Dodecahedron).

The Dynamic Geometry Software *Sketchpad v5.10* was used to construct all figures in this paper.

This work deals with the solutions of the problem of determining the minimal regular Tetrahedra which contain each of the other four classical unit polyhedra. Our results agree with the known solutions in [1] and [2], but with different methods in some cases.

**Tetrahedra Enclosing Cubes**

We begin with the unit Cube *C*, and a Tetrahedron **T** = *ABCD* which is in the corner of the first octant of the vector space $R^3$ with sides on the positive axes of length *3*. Then *A = (3,0,0), B = (0,3,0), C = (0,0,3),* and *D = (0,0,0),* Fig. 1.
Construct the unit cube **C** = *DE...K* with the following coordinates: *D = O, E = (1,0,0), F = (1,1,0), G = (0,1,0), H = (0,1,1), I = (0,0,1), J = (1,0,1), K = (1,1,1),* Fig. 2.
The face *ABC* of **T** is in the plane with equation $x + y + z = 3$, so vertex *K* is in the face *ABC*.

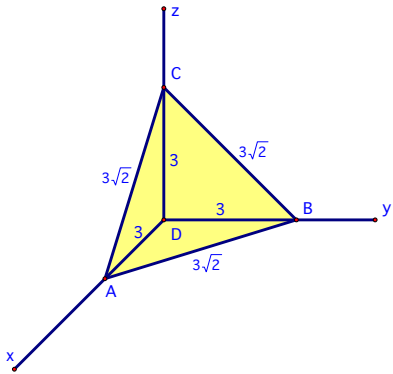
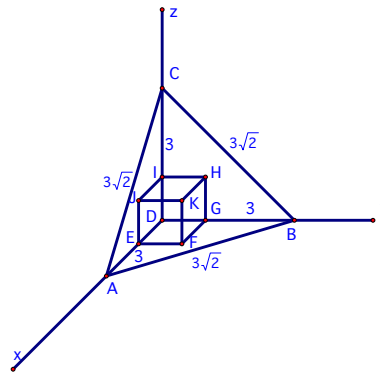

        *Figure 1*                                    *Figure 2*

The cross section at *z = 1* is a right triangle with sides of length *2*, and hypotenuse of length *2√2*, Fig. 3. It is known that this triangle has minimum area about the square *HIJK* [3], [4].

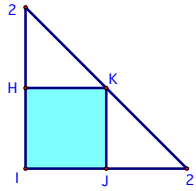
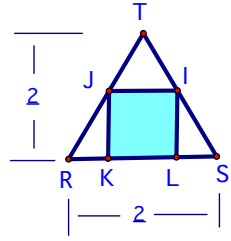

*Figure 3*  *Figure 4*

The volume of **T** = ABCD is (ABCD) = hA/3, h the height of **T**, and A the area the base, which is ABD in this case. Thus (ABCD) = 4.5 units³.

We call **T** = ABCD a *Right Tetrahedron*.

The right tetrahedron above is not regular, the sides are of *2* different lengths. If we choose *ABCD* to be a regular tetrahedron **T** with sides *s*, then the position of the unit cube **C** = *EF...KL* inside the tetrahedron will be different. So we can assume that the cube **C** is setting on the bottom side *ABC* of *ABCD*. This is called *Standard Position.*

In order to achieve a minimum size tetrahedron, we will also assume that one side, *KL,* of the top of **C**, is on side *ABD* of **T**, Fig.4.

We are also assuming that this is all happening in **R**³, so the plane *z = 1* will intersect **T** in a regular triangle *RST,* Fig. 5.

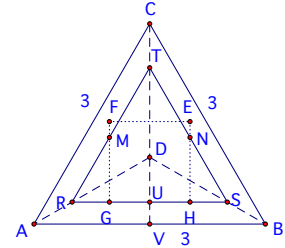

Suppose **T** has sides of length *s = 3*. Then the height is *h = (√6)s/3 = √6,* and the slope of side *ABD* to the base *ABC* is *m = (√6)/(√3/2) = 2√2*. So the distance between the base triangle *ABC* and the projection of the triangle *RST* on the *xy*-plane is *√2/4,* since the height of *RST* above the *xy*-plane is *1*. Project the point *D* onto the *xy*-plane. Since the triangle *RSD* is isosceles with base angles of *30°,* the length of *RS* is *2(√3/2-√2/4)√3 ~ 1.775.*

The base *GHEF* of **C** intersects the projected triangle *RST* at points *M* and *N*. This is because the distance from *G, H* to *M, N,* respectively, is *(RS - 1)√3/2 ~ 0.671.*

*Figure 5*

Thus the top corners of **C** do not fit inside the tetrahedron *ABCD*.

If the triangle *RST* is enlarged to contain the top square *IJKL* of **C**, the sides would have to be lengthened to equal at least *1 + 2√3/3 ~ 2.135*. This would require a corresponding lengthening of *AB* to *1+2√3/3+√6/2 ~ 3.37945*.

The volume of **T** with side *s* is *(ABCD) = (√2)s³/12 u³*, so if *s = 3*, then *(ABCD) = 9√2/4 u² ~ 3.182 u²*. If *s ~ 3.37945*, then the larger tetrahedron **T'** = *A'B'C'D'* would have volume *(A'B'C'D') ~ 4.5485 u³*.

This is larger than the volume of the right tetrahedron above, but it is the minimum *regular* tetrahedron which contains the unit cube **C**. Thus a regular minimal tetrahedron need not be a minimal tetrahedron.

These precise results on regular polyhedra agree with those obtained in [1] and [2].

Call a triangle *ABC* an *sxs triangle* if the base is of length *s* and the height is equal to *s*.

The triangles *ABD, ACD,* and *BCD* in the right tetrahedron example above are *3x3* triangles.

Similarly, call a tetrahedron *ABCD* an *sxsxs* tetrahedron if the base *ABC* is an *sxs* triangle and the height is equal to *s*. In this case $(ABCD) = hA/3 = s^3/6\ u^3$, where the area of the base *A* is $(ABC) = s^2/2\ u^2$.

The right tetrahedron above is a *3x3x3* tetrahedron, and it's area is $(ABCD) = 9/2\ u^3$.

Let *ABCD* be a *3x3x3* tetrahedron, with base *ABC* an acute triangle, and vertex *D* above the base. Assume that a unit cube *HI...NO* is in standard position in *ABCD* with *LM* in face *ACD*. Then the cross section by the plane $z = 1$ of the tetrahedron is a *2x2* triangle, say *EFG,* Fig. 6, which contains the top face of of the cube, *LMNO*. It is known that a *2x2* triangle, such as *EFG,* has minimum area about the unit square *LMNO*, and the area satisfies $(EFG) = 2\ u^2$ [3][4].

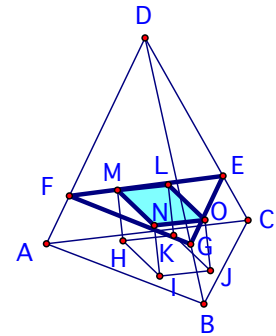

*Figure 6*

Suppose *RSTU* is another tetrahedron of minimum volume which contains the unit cube in standard position, and $(RSTU) < 9/2\ u^3$, for $h = 3$. Then the base *RST* satisfies $(RST) < 9/2\ u^2 = (ABC)$, so the cross section *VWX* at level $z = 1$ satisfies $(VWX) < (EFG) = 2\ u^2$. This contradicts the result concerning minimum area triangles which contain a unit square [3], [4]. Hence, we must have $(RSTU) \geq 9/2\ u^3$, if the other conditions are met.

Let *R'S'T'U'* be yet another tetrahedron of minimum volume which contains the unit cube in standard position, and $h \neq 3$. Then if *V'W'X'* is the cross section at level $z = 1$, the height *h'* of tetrahedron *V'W'X'U'* is not 2, $h' \neq 2$. But *V'W'X'* must be a *2x2* triangle in order to be minimal. This is impossible, for then $(V'W'X'U') > 4/3\ u^3$, so $(R'S'T'U') > 9/2\ u^3$.

In fact the minimal tetrahedron with cube in standard position, for any height $h > 1$, must always have a cross section at level $z = 1$ which is a *2x2* triangle. The base triangle will then be an *axa* triangle for some $a > 0$, by equal angles congruence.

**Tetrahedra Enclosing Octahedra**

The regular octahedron *O* with sides *s* can be thought of as gluing 2 square base regular pyramids with sides *s* together by their bases, Fig. 7. All faces are equilateral triangles, and each vertex has *4* faces meeting there. It has *8* faces and *6* vertices.

By adding regular tetrahedra with side *s* to alternate faces of the octahedron *O,* a tetrahedron *T* with sides of length *2s* is determined, Fig. 8.

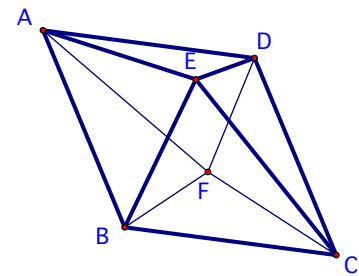

*Figure 7*

Since the sides of the octahedron *O* make up part of the sides of the tetrahedron, this is clearly the smallest regular tetrahedron which contains *O*.

If ***O*** is a unit regular octahedron, the volume of this tetrahedron is *(ABCD) = (2√2)/3 u³ ~ 0.94281 u³*, since *s = 2*. The volume of ***O*** is *√2/3 ~ 0.47140*.

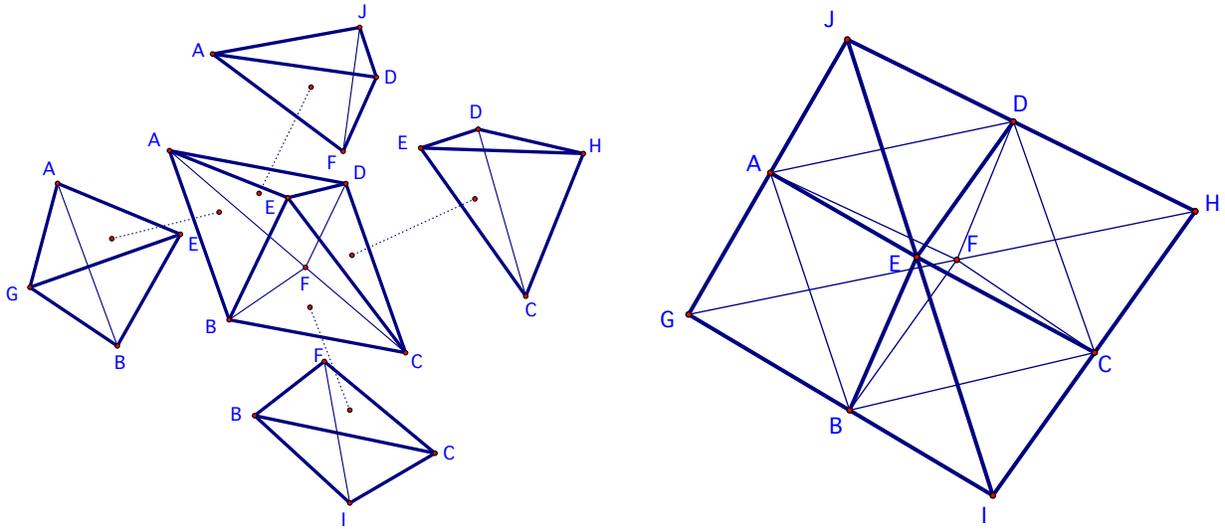

*Figure 8*

**Tetrahedra Enclosing Icosahedra**

A regular Icosahedron ***I*** with sides *s* is made up of *20* equilateral triangles such that each of the *12* vertices is the center of a pentagon of triangles, Fig. 9.
The major diagonal of *I* has length *s√(ø√5) = (s/2)√(10+2√5)*, where *ø = (1+√5)/2,* the golden ratio. This diagonal is also the diameter of the circumsphere about ***I***. We will assume *s = 1*, so that the radius of this circumsphere is *√(10+2√5)/4 ~ 0.95106*.
Let ***T*** be a tetrahedron of side *t*. Then ***T*** contains an insphere of radius *r = t√6/12*. If we set this radius equal to the radius of the circumsphere about ***I***, we have *t√6/12 = √(10+2√5)/4*, so *t = √(15 + 3√5) ~ 4.65921*.
Thus ***I*** is enclosed in ***T***, however this is not a very tightly fitting tetrahedron about ***I***, as the volume of ***I*** is *√2(15+3√5)³/12 u³ ~ 11.91982 u³*.

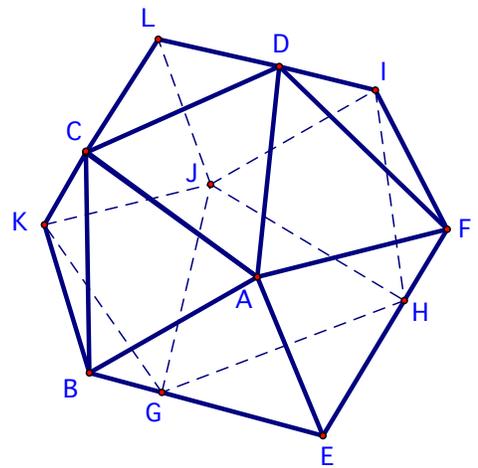

*Figure 9*

A tighter fitting tetrahedron can be constructed by building an octahedron ***O*** about ***I*** as follows. Let a *Golden Rectangle* be one which has dimensions with the ratio *1:ø*. Using *3* golden rectangles with dimensions *2x2ø* arrange them so that they are mutually perpendicular to each other, and no major axes coincide, Fig. 10. If the corners of the rectangles are connected to each other without going through the rectangles a regular icosahedron results which has sides of length *2*, Fig. 10.
The *xyz*-coordinates of the corners of the rectangles centered at the origin are: *(0, ±1, ±ø), (±ø, 0, ±1), (±1, ±ø, 0)*.

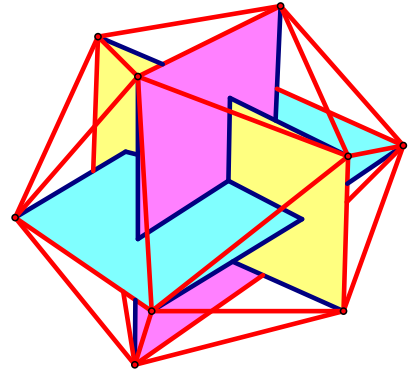

Enclosing each rectangle in a square, as shown in Fig. 11, determines an octahedron ***O*** (the corners of the squares are the vertices of ***O***) which contains the icosahedron ***I***, Fig.12, 13. The squares are *ABCD, AGCH,* and *BGDH*. Since the rectangles are of sides *2x2ø*, the squares (and octahedron) are of side length *√2(1+ø) ~ 3.7025.*
Then, we know from our studies above that an octahedron has a very close relation to a tetrahedron, so the icosahedron ***I*** here is thus enclosed in this tetrahedron ***T***.

*Figure 10*

The enclosing tetrahedron ***T*** has sides of length twice that of ***O***, so it has sides of length *2√2(1+ø) ~ 7.4049*. But the sides for ***I*** are twice those above, so this tetrahedron would have sides of length *3.7025* in comparison.
The volume of the minimum tetrahedron ***T*** then is *√2(3.7025)³/12 u³ ~ 5.9816 u³*.

This is a major reduction in size of the enclosing tetrahedron of the icosahedron, and it agrees with the minimum solutions found by Croft [1] and Firsching [2].

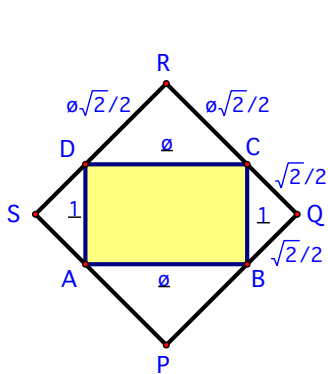 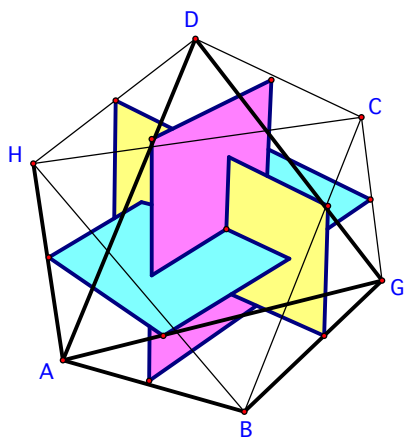 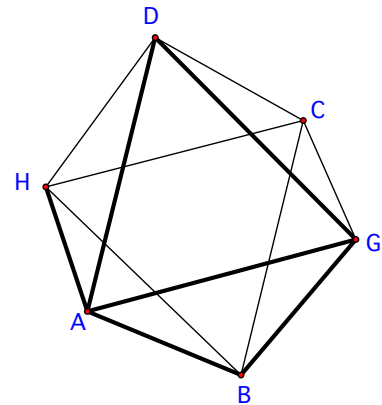

*Figure 11*         *Figure 12*         *Figure 13*

**Duals of Polyhedra**

Two polyhedra are *duals* of each other if the number of faces of one equals the number of vertices of the other, and the number of vertices of this first polyhedron also equals the number of faces of the other, the number of sides remaining the same for each polyhedron.

For example, the tetrahedron is self dual, since it has *4* vertices, *6* sides, and *4* faces.
A cube and a octahedron are duals since the cube has *8* vertices, *12* sides, and *6* faces, while the octahedron has *6* vertices, *12* sides, and *8* faces.
The most interesting case is the dodecahedron and it's dual the icosahedron.
The icosahedron has *12* vertices, *30* sides, and *20* faces. The dodecahedron has *20* vertices, *30* sides, and *12* faces.
What's more if you place a point in the center of each of the faces of an icosahedron, and connect vertices whenever the faces share a side, you get a dodecahedron. So the icosahedron contains a dodecahedron.
On the other hand, if you place a point on the center of each face on a dodecahedron you have the coordinates of *3* golden triangles which are mutually perpendicular, Fig.12, and thus determine an icosahedron. So the dodecahedron contains an icosahedron.

**Tetrahedra Enclosing Dodecahedra**

An easy way to contain a dodecahedron *D* in a tetrahedron *T* which is a relatively close fit is to form the icosahedron *I* about *D*, by duality, with vertices of *D* on the face centers of *I*, Fig. 14, and then to build the Octahedron *O* about *I*, as done above. The tetrahedron *T* which follows, with *O* an integral part of it, is a tetrahedron about *D* which is a fairly close fit. If *D* is a unit dodecahedron, then the side length of *I* is *3(√5 - 1)/2 ~ 1.8541*. The radius *r* of the circumsphere about *I* is then *r = 3(√5-1)√(10+2√5)/8 ~ 1.76336*.
If *t* is the side length of *T*, then *t = (2√6)r ~ 8.63864*.
The volume of this tetrahedron *T* is approximately *75.9749 u³*.
But it is not a very minimal tetrahedron, since the volume of the unit dodecahedron is approximately *7.6631 u³*.

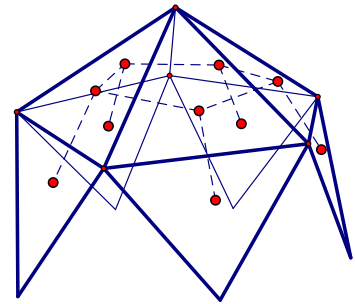

*Figure 14*

A smaller tetrahedron *T* can be obtained by assuming the icosahedron *I* is on point, that is a major axis from vertex to opposite vertex is perpendicular to the plane. Also, the enclosed unit dodecahedron *D* has a pair of horizontal opposite parallel faces, Fig.s 14, 15. Then instead of choosing the vertices of *D* in the center of the faces of *I*, we rotate *D* about the vertical axis *36°* one way or the other. This makes the vertices on the top and bottom levels move to the vertical sides of the top and bottom pentagons in *I*, and the vertices in between these top and bottom faces over to the next face, Fig. 15.

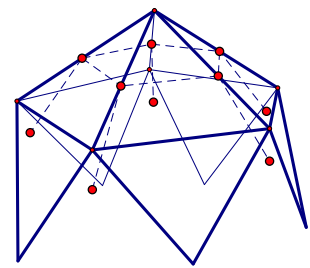

*Figure 15*

Then the length of the edge of the surrounding icosahedron *I* is equal to *(15+√5)/10 ~ 1.7236 u.*
The rest of the construction follows as given above.
The volume of *T* using this construction is equal to approximately *27.39802 u³* [2].